\documentclass[a4paper,12pt]{amsart}
\usepackage{ifthen}
\usepackage{graphicx}
\usepackage{mathrsfs}
\usepackage{color}
\nonstopmode
\numberwithin{equation}{section}
\usepackage{pgfplots}
\usepackage{tikz}
\pgfplotsset{compat=1.16}

\newtheorem{thm}{Theorem}[section]
\newtheorem{cor}[thm]{Corollary}
\newtheorem{lem}[thm]{Lemma}
\newtheorem{prop}[thm]{Proposition}
\newtheorem{rem}{Remark}

\theoremstyle{definition}
\newtheorem{defn}{Definition}

\newenvironment{pf}[1][]{%
 \vskip 3mm
 \noindent
 \ifthenelse{\equal{#1}{}}%
  {{\slshape Proof. }}%
  {{\slshape #1.} }%
 }%
{\qed\bigskip}

\newcounter{alphabet}

\newcommand{\A}{{\mathcal A}}

\newcommand{\C}{{\mathbb C}}
\newcommand{\D}{{\mathbb D}}

\newcommand{\R}{{\mathbb R}}

\newcommand{\sh}{{\mathcal S}_\mathrm{H}}
\newcommand{\T}{{\mathbb T}}

\newcommand{\spl}{{\mathcal{SP}}}

\newcommand{\harm}{{\mathrm{H}}}

\newcommand{\har}{{\mathcal{H}}}

\newcommand{\st}{{\mathcal{SS}}}

\renewcommand{\Im}{{\,\operatorname{Im}\,}}

\renewcommand{\Re}{{\,\operatorname{Re}\,}}

\newcommand{\loc}{{\operatorname{loc}}}
\renewcommand{\mod}{{\operatorname{mod}\,}}
\newcommand{\inv}{^{-1}}

\newcommand{\aand}{{\quad\text{and}\quad}}

\newcounter{minutes}
\divide\time by 60
\newcounter{hours}
\multiply\time by 60 \addtocounter{minutes}{-\time}

\begin{document}
\bibliographystyle{amsplain}
\title[Harmonic spirallike and harmonic strongly starlike functions]
{Harmonic spirallike functions and harmonic strongly starlike functions}
\def\thefootnote{}
\footnotetext{
\texttt{\tiny File:~\jobname .tex,
          printed: \number\year-\number\month-\number\day,
          \thehours.\ifnum\theminutes<10{0}\fi\theminutes}
}
\makeatletter\def\thefootnote{\@arabic\c@footnote}\makeatother
\author[X.-S. Ma]{Xiu-Shuang Ma}
\address{Graduate School of Information Sciences \\
Tohoku University\\
Aoba-ku, Sendai 980-8579, Japan}
\email{ma.xiushuang.t5@dc.tohoku.ac.jp \\
maxiushuang@gmail.com}
\author[S. Ponnusamy]{Saminathan Ponnusamy}
\address{S. Ponnusamy, Department of Mathematics \\
Indian Institute of Technology Madras \\
Chennai-600 036, India.}
\email{samy@iitm.ac.in}
\author[T. Sugawa]{Toshiyuki Sugawa}
\address{Graduate School of Information Sciences \\
Tohoku University\\
Aoba-ku, Sendai 980-8579, Japan}
\email{sugawa@math.is.tohoku.ac.jp}
\begin{abstract}
Harmonic functions are natural generalizations of conformal mappings. In recent years, a lot of work have been done by some researchers who focus on harmonic starlike functions. In this paper, we aim to introduce two classes of harmonic univalent functions of the unit disk, called hereditarily $\lambda$-spirallike functions and hereditarily strongly starlike functions, which are the generalizations of $\lambda$-spirallike functions and strongly starlike functions, respectively. We note that a relation can be obtained between this two classes. We also investigate analytic characterization of hereditarily spirallike functions and uniform boundedness of hereditarily strongly starlike functions. Some coefficient conditions are given for hereditary strong starlikeness and hereditary spirallikeness. As a simple application, we consider a special form of harmonic functions.
\end{abstract}

\thanks{
The present study was supported in part by Graduate Program in Data Science, Tohoku University.
}

\maketitle

\section{Introduction}

Logarithmic spirals frequently appear in Complex Analysis.
They are invariant under similarities, which constitutes the main reason why
they appear quite naturally.
We make a more specific definition of logarithmic spirals.
Let $\lambda$ be a real number with $|\lambda|<\pi/2.$
A plane curve of the form $w=w_0\exp(t e^{i\lambda}),~ t\in\R,$ for some \
$w_0\in \C\setminus\{0\}$ is called a \emph{$\lambda$-spiral} (about the origin).
We denote by $[0,w_0]_\lambda$ the $\lambda$-spiral segment
$\{w_0\exp(t e^{i\lambda}): t\le 0\}\cup\{0\}.$
A domain $\Omega$ in $\C$ is called $\lambda$-spirallike (with respect to the origin)
if $[0,w]_\lambda\in \Omega$
for all $w\in \Omega.$
Note that $[0, w_0]_0$ is nothing but the segment $[0, w_0]$ and thus
$0$-spirallike means \emph{starlike}.

Let $\A$ denote the class of analytic functions on the unit disk $\D=\{z\in\C: |z|<1\}$
and $\A_0, ~\A_1$ be its subclasses consisting of functions $g$ and $h$ with
$g(0)=0,~ h(0)=0, h'(0)=1,$ respectively.
A function $f$ in $\A_1$ is called $\lambda$-spirallike if $f$ maps $\D$
univalently onto a $\lambda$-spirallike domain.
It is well known (see \cite{Duren:univ}) that a function $f\in\A_1$ is
$\lambda$-spirallike if and only if the following inequality holds:
\begin{equation*}
\Re\left(
e^{-i\lambda}\frac{zf'(z)}{f(z)}
\right)>0,\quad
z\in\D.
\end{equation*}
In particular, we observe that $\lambda$-spirallikeness is a hereditary property.
Precisely speaking, if $f\in\A_1$ is $\lambda$-spirallike, then
$f_r(z)=f(rz)/r$ is again $\lambda$-spirallike for each $0<r<1.$
We denote by $\spl(\lambda)$ the class of $\lambda$-spirallike functions in $\A_1.$ See \cite{PW-14} for certain aspects of recent study on spirallike functions.

We next introduce the notion of strong starlikeness.
Let $\alpha$ be a real number with $0<\alpha<1$
and put $\tau=\tan(\pi\alpha/2).$
Let $V_\alpha$ be the Jordan domain bounded by the two logarithmic spiral segments
$\{e^{(-\tau+i)\theta}: 0\le \theta\le \pi\}$ and
$\{e^{(\tau+i)\theta}: -\pi\le \theta\le 0\}.$
Set $w_0 V_\alpha=\{w_0w: w\in V_\alpha\}.$
Note that $V_\alpha$ contains the disk $|w|<e^{-\pi\tau};$ namely,
\begin{equation}\label{eq:V}
V_\alpha\supset \left\{w: |w|<\exp\big(-\pi\tan(\pi\alpha/2)\big)\right\}.
\end{equation}
We remark that $w_0V_\alpha$ shrinks to the segment $[0, w_0)$ as $\alpha\to1.$
A domain $\Omega$ in $\C$ is called \emph{strongly starlike of order $\alpha$}
(with respect to the origin) if $w_0 V_\alpha\subset\Omega$ for all $w_0\in\Omega.$
A function $f$ in $\A_1$ is called strongly starlike of order $\alpha$
if $f$ maps $\D$ univalently onto a strongly starlike domain of order $\alpha.$
It is known (see \cite{SugawaDual}) that $f\in\A_1$ is strongly starlike of order $\alpha$
if and only if
\begin{equation}\label{eq:ss}
\left|
\arg\frac{zf'(z)}{f(z)}
\right|<\frac{\pi\alpha}2,\quad z\in\D.
\end{equation}
Thus, we also see that strong starlikeness of order $\alpha$ is a hereditary property.
We denote by $\st(\alpha)$ the class of strongly starlike functions
in $\A_1$ of order $\alpha.$
In view of the above condition, it is immediate to obtain a relation with spirallike functions
as in
\begin{equation}\label{eq:ss-spl}
\st(\alpha)=\spl\big(\tfrac{\pi(1-\alpha)}2\big)\cap
\spl\big(-\tfrac{\pi(1-\alpha)}2\big).
\end{equation}
In particular, we observe that a domain $\Omega$ is strongly starlike of order $\alpha$
if and only if $\Omega$ is $\pi(1-\alpha)/2$-spirallike and $-\pi(1-\alpha)/2$-spirallike
simultaneously.
Originally, the notion of strongly starlike functions was introduced by
Stankiewicz \cite{Stank66} and Brannan and Kirwan \cite{BK69}, independently,
with \eqref{eq:ss} being the definition.
More characterizations of strongly starlike functions are summarized in \cite{SugawaDual}.

Let $\har$ denote the class of (complex-valued) harmonic functions on $\D.$
We denote by $\har_0$ the functions $f\in\har$ normalized
by $f(0)=0$ and $f_z(0)=1.$
We note that every function $f$ in $\har_0$ can be expressed as $f(z)=h(z)+\overline{g(z)}$
for some $g\in\A_0$ and $h\in\A_1.$
According to Clunie and Sheil-Small \cite{CS84}, we denote by $\sh$ the set
of orientation-preserving harmonic univalent functions $f$ in $\har_0.$
Here, orientation-preserving means that the Jacobian $J_f=|f_z|^2-|f_{\bar z}|^2$
is positive on $\D.$
We want to extend the notions of $\lambda$-spirallike functions and strongly starlike functions
of order $\alpha$ to harmonic functions.
Before doing it, let us make a few remarks.
If $f\in\sh$ maps $\D$ onto a starlike domain with respect to the origin,
the image $f(\D_r)$ need not be starlike for some $0<r<1,$ where $\D_r=\{z: |z|<r\}.$
Indeed, we consider the harmonic Koebe function $k(z)=h(z)+\overline{g(z)}$ with
\begin{equation*}
h(z)=\frac{z-\frac12 z^2+\frac16 z^3}{(1-z)^3},\quad g(z)=\frac{\frac12 z^2+\frac16 z^3}{(1-z)^3}.
\end{equation*}
It is known (see \cite[Example 5.4]{CS84} or \cite[\S 5.3]{Duren:harm}) that $k$
maps $\D$ univalently onto the slit domain $\C\setminus (-\infty,-1/6],$ which is starlike with respect to the origin.
On the other hand, as we will see in the next section, 
$k(\D_{r_0})$ is not starlike with respect to the origin for $r_0=\sqrt{5}/3.$
Thus, starlikeness is not a hereditary property for harmonic univalent functions in $\sh.$
Therefore the notion \emph{full starlikeness} introduced by Chuaqui, Duren and Osgood 
\cite[p.~138]{CDO04} makes sense.
Here, a harmonic function $f\in\har_0$ 
is called \emph{fully starlike} if $f$ maps each circle $|z|=r~ (0<r<1)$
injectively onto a starlike curve with respect to the origin.
Note that we do not assume $f$ to be univalent on $\D.$
Indeed, they gave in \cite{CDO04} a fully starlike harmonic function $f$ 
which is not locally univalent on $\D.$
Therefore, this notion is not very convenient for our aim.
We will call $f\in\har_0$ \emph{hereditarily starlike} if $f$ is
a fully starlike harmonic univalent function on $\D.$
Note that a fully starlike harmonic function $f$ on $\D$ with non-vanishing Jacobian is (globally) univalent (see Lemma \ref{lem:AM} below).
It is natural to extend this notion to the spirallike and strongly starlike cases.

\begin{defn}
Let $\lambda$ and $\alpha$ be real numbers with $|\lambda|<\pi/2$ and $0<\alpha<1.$
A harmonic function $f$ in $\har_0$ is called \emph{hereditarily $\lambda$-spirallike}
if $f$ is orientation-preserving and univalent on $\D$ and if $f(\D_r)$ is
$\lambda$-spirallike for each $0<r<1.$
The class of such functions will be denoted by $\spl_\harm(\lambda).$
Similarly, a harmonic function $f\in \har_0$ is called
\emph{hereditarily strongly starlike of order $\alpha$} if it is
orientation-preserving and univalent on $\D$ and if $f(\D_r)$ is a strongly starlike
domain of order $\alpha$ for each $0<r<1.$
We denote by $\st_\harm(\alpha)$ the class of such functions.
\end{defn}

In particular, the class $\spl_\harm(0)$ consists of hereditarily starlike harmonic functions
on $\D.$
We would like to point out here that these classes are not considered in the literature
though spirallike logharmonic mappings and spirallike $C^1$-functions are studied by
\cite{AH87} and \cite{AM81}, respectively.

As we saw in \eqref{eq:ss-spl}, a domain $\Omega$ with $0\in\Omega\subset\C$
is strongly starlike of order $\alpha$ if and only if $\Omega$ is
$\pm\pi(1-\alpha)/2$-spirallike at the same time.
Therefore, we have similarly

\begin{equation}\label{eq:Hss-spl}
\st_\harm(\alpha)=\spl_\harm\big(\tfrac{\pi(1-\alpha)}2\big)\cap
\spl_\harm\big(-\tfrac{\pi(1-\alpha)}2\big).
\end{equation}

\section{Analytic characterization of hereditarily spirallike functions}

For continuously differentiable functions $f\in C^1(\D),$ we define the differential operator
$D$ by
$$
Df(z)=zf_z(z)-\bar zf_{\bar z}(z),
$$
where $f_z=(f_x-if_y)/2$ and $f_{\bar z}=(f_x+if_y)/2.$ Here $f_x$ and $f_y$ are the partial derivatives of $f$ with respect to $x=\Re z$ and $y=\Im z,$ respectively.
Al-Amiri and Mocanu \cite{AM81} gave a sufficient condition of $\lambda$-spirallikeness even for
functions in $C^1(\D).$
We will show that the condition is also necessary.

\begin{lem}\label{lem:AM}
Let $\lambda$ be a real number with $|\lambda|<\pi/2.$
Suppose that a function $f\in C^1(\D)$ satisfies the conditions
that $f(z)=0$ if and only if $z=0,$ and that $J_f=|f_z|^2-|f_{\bar z}|^2>0$ on $\D.$
Then $f$ is injective on $\D$ and $f(\D_r)$ is $\lambda$-spirallike for each $0<r<1$
if and only if
\begin{equation}\label{eq:hspl}
\Re\left(
e^{-i\lambda}\frac{Df(z)}{f(z)}
\right)> 0,\quad
z\in\D\setminus\{0\}.
\end{equation}
\end{lem}


For explanations, we recall a convenient quantity.
For $w\in\C\setminus\{0\},$ we will say that the $\lambda$-argument
of $w$ is $\theta$ if $w$ lies on the $\lambda$-spiral $\gamma_{\lambda,\theta}
=\{e^{i\theta}\exp(t e^{i\lambda}): t\in\R\}.$
We will write $\arg_\lambda w=\theta$ in this case.
Note that the $\lambda$-argument is determined up to an integer multiple of $2\pi$
and a more explicit expression is available as follows:
$$
\arg_\lambda w=\arg w-(\tan\lambda)\log|w|\quad (\mod 2\pi).
$$
This terminology was introduced in \cite{KS12spl} but the same idea was essentially used
in \cite{AM81} and other papers earlier.

\begin{pf}[Proof of Lemma $\text{\ref{lem:AM}}$]
As we mentioned before, the ``if~" part was shown by Al-Amiri and Mocanu \cite{AM81}.
For completeness, we describe the essential ideas for this part.
Let $C_r=f(\partial\D_r)$ for $0<r<1.$
Note that each $C_r$ does not pass through the origin by assumption.
We will show that $\{C_r\}$ is a family of non-intersecting Jordan curves.
Since $C_r$ has winding number 1 about the origin, we may take a continuous
branch of $\phi(\theta)=\arg_\lambda f(re^{i\theta})$ with period relation
$\phi(\theta+2\pi)=\phi(\theta)+2\pi.$
A straightforward computation (see \cite[p.~63]{AM81}) leads to
\begin{equation}\label{eq:phi}
\phi'(\theta)=\frac1{\cos\lambda}
\Re\left(e^{-i\lambda}\frac{Df(z)}{f(z)}\right)>0.
\end{equation}
Hence, $\phi(\theta)$ is (strictly) increasing,
which implies that $f$ is injective on each circle $|z|=r;$
in other words, $C_r$ is a Jordan curve, and that the inside of $C_r$
is a $\lambda$-spirallike domain.

Now, we need only to show that $C_r$ lies in the Jordan domain bounded by $C_{r'}$
for $0<r<r'<1.$
To this end, fix $\phi\in\R$ and we express the unique intersection point
of $C_r$ and $\gamma_{\lambda,\phi}$ as $f(re^{i\theta})=\exp(i\phi+te^{i\lambda})$
for $t=t(r)\in\R$ and $\theta=\theta(r)\in\R.$
Then, it suffices to check that $t(r)<t(r')$ for $0<r<r'<1.$
By (12) in \cite{AM81} or by a formal computation, we obtain the relation
\begin{equation}\label{eq:J}
|f(z)|^2\frac{dt}{dr} \Re\left(e^{-i\lambda}\frac{Df(z)}{f(z)}\right)=r J_f(z),
\end{equation}
where $z=re^{i\theta}.$
Since $J_f>0$ by assumption, we conclude that $t=t(r)$ is increasing in $0<r<1.$
Thus we have shown the ``if~" part.

Secondly, we show the ``only if~" part.
Assume that $f$ is univalent on $\D$ and that $f(\D_r)$ is $\lambda$-spirallike for $0<r<1.$
Then the intersection of $C_r=\partial f(\D_r)$ with $\gamma_{\lambda,\phi}$ is connected
for each $0<r<1$ and $\phi\in\R$ so that $\phi(\theta)=\arg_\lambda f(re^{i\theta})$
is non-decreasing in $\theta.$
Also, $t=t(r)$ defined above is non-decreasing in $0<r<1$ and thus $dt/dr\ge 0.$
In view of \eqref{eq:phi} and \eqref{eq:J}, we obtain \eqref{eq:hspl} because $J_f>0$
by assumption.
\end{pf}

We restate the lemma in the case when $f$ is harmonic.

\begin{cor}\label{cor:spl}
Let $\lambda$ be a real number with $|\lambda|<\pi/2.$
Suppose that a function $f\in\har_0$ satisfies the conditions
that $f(z)\ne0$ for $0<|z|<1$ and
that $J_f=|f_z|^2-|f_{\bar z}|^2>0$ on $\D.$
Then $f\in\spl_\harm(\lambda)$ if and only if the inequality
\eqref{eq:hspl} holds.
\end{cor}

In particular, by \eqref{eq:Hss-spl},
we obtain the following characterization of hereditarily strongly starlike
functions of order $\alpha.$

\begin{cor}\label{cor:st}
Let $\alpha$ be a real number with $|\alpha|<1.$
Suppose that a function $f\in\har_0$ satisfies the conditions
that $f(z)\ne0$ for $0<|z|<1$ and
that $J_f=|f_z|^2-|f_{\bar z}|^2>0$ on $\D.$
Then $f\in\st_\harm(\alpha)$ if and only if
\begin{equation}\label{eq:hss}
\left|
\arg\frac{Df(z)}{f(z)}
\right|<\frac{\pi\alpha}2,\quad
z\in\D\setminus\{0\}.
\end{equation}
\end{cor}

\begin{pf}[Proof of non-hereditary starlikeness of $k(z)$]
We now show that the harmonic Koebe function $k(z)$ 
is not hereditarily starlike.
By virtue of Lemma \ref{lem:AM}, it is enough to check that the function
$k$ does not satisfy the condition $\Re[Dk/k]>0$ on $\D.$
Here,
$$
Dk(z)=zh'(z)-\overline{zg'(z)}
=\frac{z(1+z)}{(1-z)^4}-\frac{\bar z^2(1+\bar z)}{(1-\bar z)^4}.
$$
Let $z_0=(1+2i)/3\in\D.$
Then, straightforward computations yield $k(z_0)=(-17+9i)/24$
and $Dk(z_0)=-15(1+2i)/16.$
Hence, we see that $Dk(z_0)/k(z_0)=9(-1+43i)/148$ has
negative real part.
\end{pf}

Let $r_1$ be the radius of hereditary starlikeness for the harmonic
Koebe function $k.$
Then, numerical computations suggest that $0.572154<r_1<0.572155.$

\section{Uniform boundedness of hereditarily strongly starlike functions}

Brannan and Kirwan \cite{BK69} showed that a function $f\in\st(\alpha)
~(0<\alpha<1)$ admits the sharp estimate
$$
|f(z)|<|z|M(\alpha)<M(\alpha),\quad 0<|z|<1,
$$
where
$$
M(\alpha)
=\exp\left\{ 2\alpha\sum_{k=0}^\infty\frac1{(2k+1)(2k+1-\alpha)}\right\}
=\frac14\exp\left\{-\psi\big((1-\alpha)/2\big)-\gamma\right\},
$$
$\psi(x)=\Gamma'(x)/\Gamma(x)$ is the digamma function and
$\gamma=0.5772\dots$ is the Euler-Mascheroni constant.
In this section, we extend it to hereditarily strongly starlike harmonic functions of
order $\alpha.$
To this end, we recall the following result due to Hall \cite{Hall85}
(see also \cite[\S 6.2]{Duren:harm}).

\begin{lem}\label{lem:Hall}
Let $f\in\sh.$
Then there is a point $w_0\in\C$ with $|w_0|\le \pi/2$ such that
$w_0\notin f(\D).$
The bound $\pi/2$ is sharp.
\end{lem}

Our result in this section is the following.

\begin{thm}\label{thm:bdd}
Let $\alpha$ be a real number with $0<\alpha<1.$
For each $f\in\st_\harm(\alpha),$ the inequality
$
|f(z)|\le N(\alpha), z\in\D,
$
holds, where
$$
N(\alpha)=\frac\pi 2\exp\big\{\pi\tan(\pi\alpha/2)\big\}.
$$
\end{thm}

\begin{pf}
We define $f_r$ by $f_r(z)=f(rz)/r.$
Then $f_r\in \st_\harm(\alpha)\subset \sh$ for each $0<r<1.$
Let $\Omega_r=f_r(\D)$ for $0<r<1.$
Then $\Omega_r$ is a strongly starlike domain of order $\alpha.$
For an arbitrary point $w\in\Omega_r\setminus\{0\},$ we have
$w V_\alpha\subset \Omega_r.$
On the other hand, by Lemma \ref{lem:Hall}, there is a point
$w_0\in \C\setminus\Omega_r$ with $|w_0|\le \pi/2.$
In view of the relation \eqref{eq:V}, we have
$$
|w|\exp\big(-\pi\tan(\pi\alpha/2)\big)\le |w_0|\le \frac\pi 2
$$
for $w\in \Omega_r,$ which implies $|w|\le N(\alpha).$
Since $0<r<1$ was arbitrary, we have the expected conclusion.
\end{pf}

\begin{figure}
\begin{center}
\includegraphics[width=.65\textwidth]{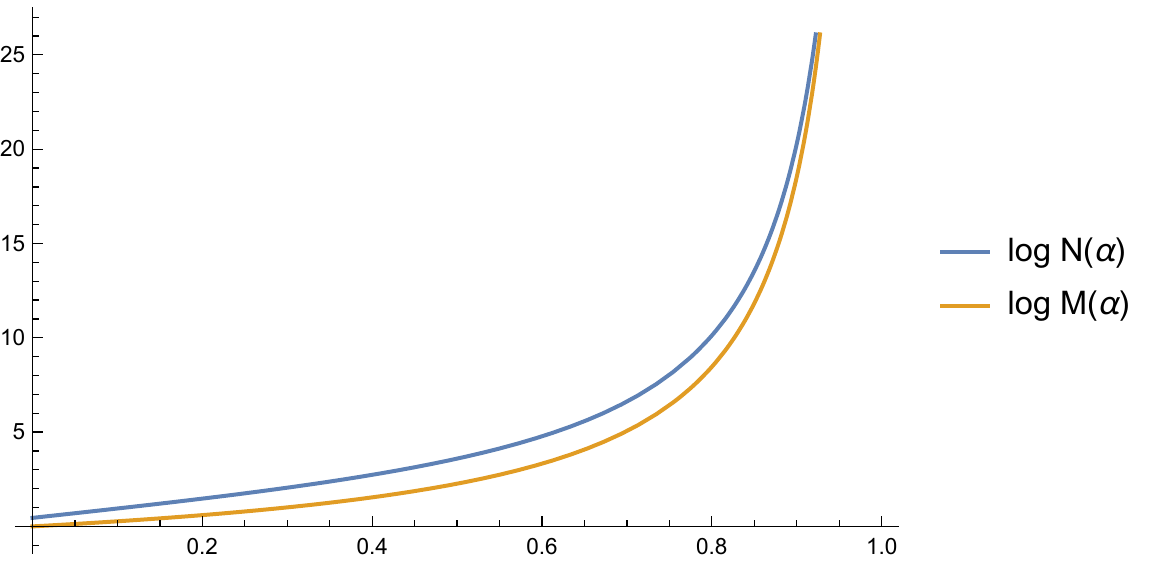}
\caption{The graph of $\log M(\alpha)$ and $\log N(\alpha)$}\label{fig:1}
\end{center}
\end{figure}

We exhibit the graph of $\log M(\alpha)$ and $\log N(\alpha)$ in Figure \ref{fig:1}.
Though $M(\alpha), N(\alpha)\to+\infty$ as $\alpha\to1,$ the graph
suggests that $\log N(\alpha)-\log M(\alpha)$ is bounded.
Indeed, that is true.
Consider the ratio 
$$
\frac{N(\alpha)}{M(\alpha)}
=2\pi\exp\big\{\pi\tan(\pi\alpha/2)+\psi((1-\alpha)/2)+\gamma\big\}
=2\pi\exp\big\{\pi\cot(\pi t)+\psi(t)+\gamma\big\},
$$
where $t=(1-\alpha)/2.$
Since $\cot x=1/x+O(x)$ and $\psi(x)=1/x-\gamma+O(x)$ as $x\to 0,$
we have $\pi\cot(\pi t)+\psi(t)+\gamma=O(t)$ as $t\to 0.$
Hence,
$$
\lim_{\alpha\to 1}\frac{N(\alpha)}{M(\alpha)}=2\pi.
$$
By numerical computations, we observed that $N(\alpha)\le 2\pi M(\alpha)$
for $0<\alpha<1.$


As an application of the boundedness, we establish quasiconformal extendability
of hereditarily strongly starlike harmonic functions under a mild condition.

First, we recall that a homeomorphism $f:\Omega\to\Omega'$ between plane domains is called $K$-quasiconformal if $f$ belongs to the Sobolev class $W_\loc^{1,2}(\Omega)$ and if the inequality $|f_{\bar z}|\le k |f_z|$ holds a.e.~on $\Omega$, where $k=(K-1)/(K+1)\in[0,1)$.
When $\Omega=\Omega'$, we call $f$ a $K$-quasiconformal endomorphism of $\Omega$.
It is well known \cite{Ahlfors:qc2} that $f_1\circ f_2$ is $K_1K_2$-quasiconformal whenever
$f_j$ is $K_j$-quasiconformal for $j=1$, $2$. A bounded domain $\Omega$ is called a \emph{$K$-quasidisk} if $\Omega=f(\D)$ for a $K$-quasiconformal mapping $f:\C\to\C$. Fait, Krzy\. z and Zygmunt \cite{FKZ76} showed the following.

\begin{lem}\label{lem:FKZ}
Let $0<\alpha<1.$
A strongly starlike function in $\st(\alpha)$ extends to a 
$\cot^2\frac{\pi(1-\alpha)}4$-quasiconformal endomorphism of $\C.$
In particular, a strongly starlike domain of order $\alpha$ is a 
$\cot^2\frac{\pi(1-\alpha)}4$-quasidisk.
\end{lem}

We extend this result to the class $\st_\harm(\alpha)$ of hereditarily strongly starlike harmonic functions of order $\alpha$.

\begin{thm}
Let $f=h+\bar g\in \st_\harm(\alpha)$ for some $0<\alpha<1.$
Suppose that the second complex dilatation $\omega=g'/h'$ of $f$ satisfies the inequality
$|\omega|\le (K-1)/(K+1)$ on $\D$ for a constant $K\ge 1.$
Then $f$ extends to a $K\cot^2\frac{\pi(1-\alpha)}4$-quasiconformal endomorphism
of $\C.$
\end{thm}

\begin{pf}
Let $\Omega=f(\D).$
By definition, $\Omega$ is a strongly starlike domain of order $\alpha.$
Let $\mu=f_{\bar z}/f_z=\overline{g'}/h'$ be the complex dilatation of $f.$
Then $|\mu|=|\omega|\le (K-1)/(K+1)<1.$
Let $w:\D\to\D$ be a quasiconformal homeomorphism with $w(0)=0, w(1)=1$ and
$w_{\bar z}/w_z=\mu$ a.e.~on $\D.$
Note that existence of such a mapping is guaranteed by the measurable Riemann
mapping theorem (see \cite{Ahlfors:qc2}).
Moreover, the mapping $w$ extends to a $K$-quasiconformal mapping of $\C$
with the property $1/w(1/z)=w(z)$ for $z\in\D.$
Then the composed mapping $F=f\circ w\inv:\D\to \Omega$ is analytic and
satisfies $F(0)=0.$
Let $a=F'(0)$ and $G=F/a.$
Since the image $G(\D)=\Omega/a$ is strongly starlike of order $\alpha,$
we observe that $G\in\st(\alpha).$
Now Lemma \ref{lem:FKZ} implies that $G$ extends to a 
$\cot^2\frac{\pi(1-\alpha)}4$-quasiconformal endomorphism of $\C.$
Hence, $f=F\circ w$ extends to a 
$K\cot^2\frac{\pi(1-\alpha)}4$-quasiconformal endomorphism of $\C$
as required.
\end{pf}

\section{Coefficient conditions for hereditary strong starlikeness}

Since the quotient $Df(z)/f(z)$ is not necessarily harmonic, it is not easy to check
conditions \eqref{eq:hspl} and \eqref{eq:hss} for a specific function $f$ in $\har_0.$
In this section, we give simple sufficient conditions in terms of the coefficients
of $f$ by employing the ideas due to Silverman \cite{Sil98}.

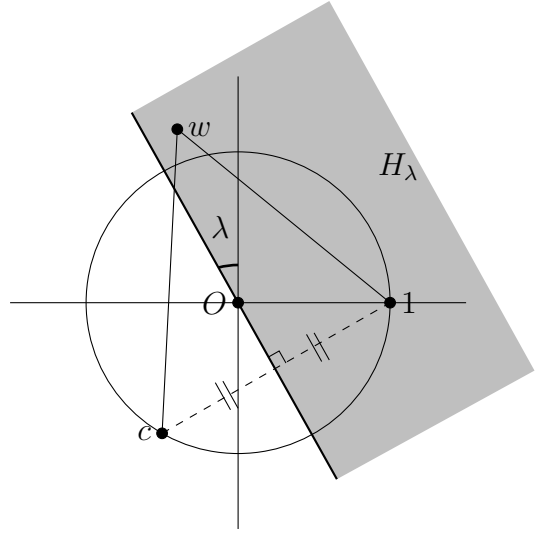
\begin{figure}
\begin{center}
\begin{tikzpicture}
\fill[lightgray] (1.2,4.0)--(-1.4,2.52)--(1.3,-2.34)--(3.9,-0.9);
\draw[thin] (-3,0) -- (3,0);
\draw[thin] (0,-3) -- (0,3);
\draw[black,line width=1] (0,0.5) arc (90:105:1);
\draw[thick] (-1.4,2.52) -- (1.3,-2.34);
\draw[thin] (0,0) circle (2cm);
\coordinate (o) at(0,0) node[left] at (o) {$O$};
\coordinate (a) at(2,0) node[right] at (a) {$1$};
\coordinate (c) at(-1,-1.732) node[left] at (c) {$c$};
\coordinate (w) at(-0.8,2.3) node[right] at (w) {$w$};
\coordinate (l) at(-0.5,1) node[right] at (l) {$\lambda$};
\coordinate (h) at(1.7,1.8) node[right] at (h) {$H_\lambda$};
\coordinate (x) at(0.4,-0.72);
\coordinate (y) at(0.54,-0.64);
\coordinate (z) at(0.63,-0.80);
\draw[thin] (a) -- (w);
\draw[thin] (c) -- (w);
\draw (x) -- (y) -- (z);
\draw[dashed] (c) -- (a);
\draw (-0.3,-1.06) -- (-0.1,-1.42);
\draw (-0.2,-1.04) -- (0,-1.4);
\draw (1.1,-0.78) -- (0.9,-0.42);
\draw (1.2,-0.76) -- (1.0,-0.40);
\filldraw [black] (0,0) circle (2pt)(2,0) circle (2pt) (-1,-1.732) circle (2pt) (-0.8,2.3)circle (2pt);
\end{tikzpicture}
\end{center}
\caption{The half-plane $H_\lambda$ and the point $c$}
\label{fig:2}
\end{figure}

We first observe that the condition \eqref{eq:hspl} means that
the quantity $Df(z)/f(z)$ lies in the half-plane $H_\lambda=
\{w\in\C: \Re(e^{-i\lambda} w)>0\}.$
Let $c$ be the mirror image of the point $1$ in the line $\partial H_\lambda.$
More precisely,
\begin{equation*}
c=-e^{2i\lambda}=-\cos{2\lambda}-i\sin{2\lambda}.
\end{equation*}
Then a point $w\in\C$ lies in the half-plane $H_\lambda$ if and only if
$|w-1|<|w-c|.$
See Figure \ref{fig:2}.
We apply this idea to deduce our result.

For $0<\alpha<1,$ we introduce the following quantities for integers $n\ge 1:$
\begin{align*}
A_n(\alpha)&=n-1+|n-e^{-i\pi\alpha}|=
n-1+\sqrt{n^2-2n\cos{\pi\alpha}+1}, \\
B_n(\alpha)&=n+1+|n+e^{i\pi\alpha}|
=n+1+\sqrt{n^2+2n\cos{\pi\alpha}+1}.
\end{align*}

\begin{lem}\label{lem:starlike}
For $n\ge 2,$ the following inequalities hold:
\begin{equation}
2n\sin\frac{\pi\alpha}2< A_n(\alpha)<B_n(\alpha)
\quad (0<\alpha<1).
\label{eq:ab}
\end{equation}
\end{lem}

\begin{pf}
Let $a=\sin{(\pi\alpha/2)}.$ Then $0<a<1$. The inequalities \eqref{eq:ab} can be written in the form
\begin{equation*}
2na< (n-1)+\sqrt{(n-1)^2+4na^2}<(n+1)+\sqrt{(n+1)^2-4na^2}.
\end{equation*}
By the triangle inequality $|n-e^{-i\pi\alpha}|-|n+e^{i\pi\alpha}|<2,$
the second inequality in \eqref{eq:ab} can be checked easily.
So we only prove the first inequality.
Indeed, we can show the stronger inequality
$$
|2na-(n-1)|  <\sqrt{(n-1)^2+4na^2},
$$
which is equivalent to
$$
\big[(n-1)^2+4na^2\big]-|2na-(n-1)|^2
=4na(n-1)(1-a)>0.
$$
Now we can check easily the first one for $n\geq 2$ and $a<1$.
\end{pf}

We are now in a position to state our main result in this section.

\begin{thm}\label{thm:coef}
Let $f=h+\bar g\in\har_0$ for $h(z)=z+a_2 z^2+a_3 z^3+\cdots$
and $g(z)=b_1z+b_2z^2+b_3z^3+\cdots.$
Suppose that the inequality
\begin{equation}\label{eq:ci}
\sum_{n=2}^{\infty}A_n(\alpha)|a_n|
+\sum_{n=1}^{\infty}B_n(\alpha)|b_n|\leq 2\sin\frac{\pi\alpha}2
\end{equation}
holds. Then $f\in\st_\harm(\alpha).$
\end{thm}

\begin{pf}
Obviously, we can assume that $f(z)$ is not identically $z.$
We first show that $f$ is orientation-preserving.
Indeed, by Lemma \ref{lem:starlike}, the condition \eqref{eq:ci} leads to
\begin{equation}\label{eq:Sil}
\sum_{n=2}^{\infty}n|a_n|
+\sum_{n=1}^{\infty}n|b_n|< 1.
\end{equation}
Therefore,
\begin{align*}
|f_z(z)|-|f_{\bar z}(z)|&=|h'(z)|-|g'(z)|
\ge 1-\sum_{n=2}^{\infty}n|a_n||z|^{n-1}
-\sum_{n=1}^{\infty}n|b_n||z|^{n-1} \\
&\ge 1-\sum_{n=2}^{\infty}n|a_n|
-\sum_{n=1}^{\infty}n|b_n|> 0
\end{align*}
for $z\in\D.$
Hence, $J_f=|f_z|^2-|f_{\bar z}|^2>0,$ 
which means that $f$ is orientation-preserving.

Let $\lambda$ be a real number with $|\lambda|<\pi/2$
and let $c=-e^{2i\lambda}$ as above.
Noting that $Df(z)/f(z)\in H_\lambda$ if and only if $|Df(z)-f(z)|<|Df(z)-cf(z)|,$
we deduce that for $0<|z|<1$

\vspace{8pt}
$|Df(z)-f(z)| < |Df(z)-cf(z)|$
\begin{align*}
 \Leftrightarrow &\left |zh'(z)-\overline{zg'(z)}-h(z)-\overline{ g(z)}\right |<
\left |zh'(z)-\overline{zg'(z)}-ch-c\overline{ g(z)}\right | \\
        \Leftrightarrow & \left|\sum_{n=2}^{\infty}(n-1)a_n z^n
-\sum_{n=1}^{\infty}(n+1)\bar b_n \bar z^n\right|
< \left|(1-c)z+\sum_{n=2}^{\infty}(n-c)a_n z^n
-\sum_{n=1}^{\infty}(n+c)\bar b_n \bar z^n\right| \\
        \Leftarrow & \left|\sum_{n=2}^{\infty}(n-1)a_n z^n\right|
+\left|\sum_{n=1}^{\infty}(n+1)\bar b_n \bar z^n\right|
< |1-c| |z|-\left|\sum_{n=2}^{\infty}(n-c)a_n z^n\right|
-\left|\sum_{n=1}^{\infty}(n+c)\bar b_n \bar z^n\right| \\
        \Leftrightarrow & \left|\sum_{n=2}^{\infty}(n-1)a_n z^n\right|
+\left|\sum_{n=1}^{\infty}(n+1)\bar b_n \bar z^n\right|
 + \left|\sum_{n=2}^{\infty}(n-c)a_n z^n\right|
+\left|\sum_{n=1}^{\infty}(n+c)\bar b_n \bar z^n\right|
         < |1-c| |z| \\
        \Leftarrow & \sum_{n=2}^{\infty}(n-1+|n-c|)|a_n|
+ \sum_{n=1}^{\infty}(n+1+|n+c|)|b_n|    \leq |1-c|=2\cos\lambda.
\end{align*}
Thus we have seen that the last inequality is sufficient for the condition
$Df/f\in H_\lambda.$
We now observe that the last inequality remains invariant when
$\lambda$ is replaced by $-\lambda.$
Therefore, taking $\lambda=\pi(1-\alpha)/2,$ we obtain the
required conclusion with the help of Corollary~\ref{cor:st}.
\end{pf}

We remark that \eqref{eq:Sil} is the condition for hereditary starlikeness
given by Silverman \cite{Sil98} (at least when $b_1=0$).
This condition also ensures that $\Re h'(z)>|g'(z)|$ in $\mathbb{D}$ 
and hence $f=h+\bar g\in\har_0$
is (hereditarily) starlike and close-to-convex with the identity function. 
See \cite[Corollary 1.4]{PYY-13} and \cite[Lemma 2.1]{KPV14}.

Let $f$ and $F$ be two harmonic functions on $\D$ of the forms
$$f(z)=\sum_{n=0}^{\infty}a_{n} z^n+\sum_{n=1}^{\infty}b_{n} \bar z^n
~~\mbox{ and }~
F(z)=\sum_{n=0}^{\infty}A_{n} z^n+\sum_{n=1}^{\infty}B_{n} \bar z^n.
$$
Then the (harmonic) convolution $f*F$ of $f$ and $F$ is defined as
\begin{equation*}
(f*F)(z)=
f(z)*F(z)=\sum_{n=0}^{\infty}a_{n}A_{n} z^n
+\sum_{n=1}^{\infty}b_{n}B_{n} \bar z^n.
\end{equation*}
In other words, for $f=h+\bar g,~ F=H+\bar G,$ the convolution is defined by
$f*F=h*H+\overline{g*G}.$
Inspired by \cite{SST78}, we have the following convolution theorem for harmonic spirallike functions.
It is useful below to note that $h(z)=h(z)*\frac{z}{1-z}$ and
$zh'(z)=h(z)*\frac{z}{(1-z)^2}$
for an analytic function $h$ on $\D$ with $h(0)=0.$

\begin{thm}\label{thm:conv}
Let $-\pi/2<\lambda<\pi/2.$
Suppose that a function $f=h+\overline{g}\in\har_0$ satisfies
$f(z)\ne0$ for $0<|z|<1$ and $J_f(z)>0$ for $z\in \D.$
Then $f\in\spl_\harm(\lambda)$ if and only if
\begin{equation}\label{eq:conv1}
(f*\varphi_{\lambda,\zeta})(z)\ne 0
\quad \text{for}~ z\in\D\setminus\{0\},\ \zeta\in \T\setminus\{-1\},
\end{equation}
where $\T$ denotes the unit circle $\partial\D$ and
\begin{equation*}
\varphi_{\lambda, \zeta}(z)= 
\frac{(1+e^{2i\lambda})z+(\zeta-e^{2i\lambda})z^2}{(1-z)^2} 
+\frac{(-1+e^{2i\lambda}-2\zeta)\bar z+(\zeta-e^{2i\lambda})\bar z^2}{(1-\bar z)^2}.
\end{equation*}
\end{thm}


\begin{pf}
By Lemma \ref{lem:AM}, $f$ belongs to $\spl_\harm(\lambda)$
if and only if $f$ satisfies the inequality \eqref{eq:hspl}.
It may be expressed in the form
\begin{equation*}
\Re \frac{1}{\cos \lambda}\left ( e^{-i\lambda}\frac{Df(z)}{f(z)}+i\sin\lambda \right ) >0
\end{equation*}
for $z\in\D\setminus\{0\}.$
The key fact in the proof is that the M\"obius transformation 
$z\mapsto \frac{z-1}{z+1}$ maps the unit circle $\T$ onto 
the extended imaginary axis $i\R\cup\{\infty\}.$
Thus, the above condition is equivalent to
$$
 \frac{1}{\cos \lambda}\left ( e^{-i\lambda}\frac{Df(z)}{f(z)}+i\sin\lambda \right )  \neq\frac{\zeta-1}{\zeta+1},\quad \zeta\in\T\setminus\{-1\}.
$$
This may be further
rephrased as
\begin{align*}
& \quad \ [Df(z)+ie^{i\lambda}(\sin\lambda) f(z)](\zeta+1)
-e^{i\lambda}(\cos\lambda)f(z)(\zeta-1) \\
&=(\zeta+1)Df(z)+(e^{2i\lambda}-\zeta)f(z) \\
&= (\zeta+1)(zh'(z)-\overline{zg'(z)})+(e^{2i\lambda}-\zeta)(h(z)+\overline{g(z)}) \\
&= [(\zeta+1)zh'(z)+(e^{2i\lambda}-\zeta)h(z)]-
[(\zeta+1)\overline{zg'(z)}+(\zeta-e^{2i\lambda})\overline{g(z)}] \\
&= h(z)*\left[\frac{(\zeta+1)z}{(1-z)^2}+\frac{(e^{2i\lambda}-\zeta)z}{1-z}\right]
-\overline{g(z)}*\left[\frac{(\zeta+1)\bar z}{(1-\bar z)^2}+\frac{(\zeta-e^{2i\lambda})\bar z}{1-\bar z}\right] \\
&=h(z)*\frac{(1+e^{2i\lambda})z+(\zeta-e^{2i\lambda})z^2}{(1-z)^2} 
+\overline{g(z)}*\frac{(-1+e^{2i\lambda}-2\zeta)\bar z+(\zeta-e^{2i\lambda})\bar z^2}{(1-\bar z)^2} \\
&=(f*\varphi_{\lambda,\zeta})(z)\ne 0
\end{align*}
\end{pf}

Taking $\pm\pi(1-\alpha)/2$ as $\lambda$ in \eqref{eq:conv1}, 
with the help of \eqref{eq:Hss-spl}, we obtain the following result.

\begin{cor}\label{cor:2}
Let $f$ be an orientation-preserving harmonic function in $\har_0$ satisfying the condition
$f(z)\ne0$ for $0<|z|<1.$ 
For $0<\alpha<1$, $f\in\st_\har(\alpha)$ if and only if
$$
(f*\varphi_{\frac{\pi(1-\alpha)}2,\zeta})(z)\ne 0
\aand
(f*\varphi_{-\frac{\pi(1-\alpha)}2,\zeta})(z)\ne 0
$$
for all $z\in\D\setminus\{0\}$ and $\zeta\in\T\setminus\{-1\}.$
\end{cor}



As a simple application of the above results, we examine hereditary strong starlikeness
of the harmonic function $f_{b,n}$ of the special form
$$
f_{b,n}(z)=z+b\overline{z}^n
$$
for $b\in\C$ and $n=1,2,3,\dots.$

\begin{prop}\label{prop:shpbd}
Let $0<\alpha<1$ and set $\lambda=\pi(1-\alpha)/2$.
Then the following are equivalent:
\begin{itemize}
\item[(i)] $f_{b,n}\in\st_\harm(\alpha)$;\vspace{4pt}
\item[(ii)] $f_{b,n}\in\spl_\harm(\lambda)$;
\item[(iii)] $|b|\le C_n(\alpha)$, where $\displaystyle
C_n(\alpha)=\frac{2\sin(\pi\alpha/2)}{n+1+|n+e^{i\pi\alpha}|}.$
\end{itemize}
\end{prop}

\begin{rem}
Since $C_n(\alpha)=2\sin(\pi\alpha/2)/B_n(\alpha),$
Lemma \ref{lem:starlike} implies that $C_n(\alpha)<1.$
\end{rem}

\begin{pf}
(i) $\Rightarrow$ (ii). It is obvious by the relation \eqref{eq:ss-spl}.

\noindent
(ii) $\Rightarrow$ (iii).  Assume that $f_{b,n}\in\spl_\harm(\lambda)$. 
By Theorem \ref{thm:conv}, $f_{b,n}$ must satisfy the condition \eqref{eq:conv1};
namely,
$$
(f_{b,n}*\varphi_{\lambda, \zeta})(z)=
(1-e^{2i\lambda})z-b[(n+1)\zeta+n+e^{2i\lambda}]\bar{z}^n\neq 0,
$$
for $0<|z|<1$ and $|\zeta|=1$ with $\zeta\ne -1.$
This implies
$$
|1-e^{2i\lambda}|\geq|b|\,|(n+1)\zeta+n+e^{2i\lambda}|,\quad \zeta\in\T\setminus\{-1\}.
$$
Hence,
$$
|b|\le \sup_{|\zeta|=1,\zeta\ne -1}\frac{|1-e^{2i\lambda}|}{|(n+1)\zeta+n+e^{2i\lambda}|}
=C_n(\alpha).
$$
(iii) $\Rightarrow$ (i). 
Condition (iii) means the inequality
$
B_n(\alpha)|b|\le2\sin(\pi\alpha/2).$
We now conclude that $f_{b,n}\in\st_\harm(\alpha)$ by Theorem \ref{thm:coef}.
\end{pf}

\def\cprime{$'$} \def\cprime{$'$} \def\cprime{$'$}
\providecommand{\bysame}{\leavevmode\hbox to3em{\hrulefill}\thinspace}
\providecommand{\MR}{\relax\ifhmode\unskip\space\fi MR }
\providecommand{\MRhref}[2]{%
  \href{http://www.ams.org/mathscinet-getitem?mr=#1}{#2}
}
\providecommand{\href}[2]{#2}

\end{document}